\documentclass{amsart}
\usepackage{hyperref}  

\hypersetup{
pdftitle = {A determinant-like formula for the Kostka numbers},
pdfsubject = {Fakultaet fuer Mathematik, Universitaet Bielefeld, Germany},
pdfauthor = {Mathias Lederer},
pdfkeywords = {Kostka, FWF, P16641, Representations of finite symmetric groups, exact enumeration problems, partitions of integers, Lederer},
}
\include{aoc-01.sty} 
\newcommand{\Z}{{\,\mathbb{Z}\,}}
\newcommand{\N}{{\,\mathbb{N}\,}}
\newcommand{\Q}{{\,\mathbb{Q}\,}}

\newtheorem{thm}{Theorem}
\newtheorem{dfn}{Definition}
\newtheorem{lmm}{Lemma}
\newtheorem{pro}{Proposition}

\begin{document}

\title{A determinant-like formula for the Kostka numbers}

\author{Mathias Lederer}
\address{Fakult\"at f\"{u}r Mathematik, Universit\"{a}t Bielefeld, Bielefeld, Germany}
\thanks{This research was partly supported by project P16641 of FWF (Austrian Science Fund)}
\email{mlederer@mathematik.uni-bielefeld.de}


\date{January 10, 2004}


\begin{abstract}
  Young tableaux are ubiquitous in various branches of mathematics. There are two counting formulas for standard Young tableaux.
  The first involves a determinant and goes back to Frobenius and Young, and the second is the hook formula by Frame, Robinson and Thrall. 
  We present a generalization of the determinant formula for semistandard Young tableaux of given shape and of given content. Our counting formula 
  -- though not a determinant -- is a generalization of the determinant formula by Frobenius and Young. 
\end{abstract}

\maketitle

\section{Introduction}\label{secintro}

Let $\alpha=(\alpha_1,\alpha_2,\ldots)$ be a proper partition of the integer $n\geq1$, 
i.e. a sequence of integers $\alpha_i\geq0$ such that $\alpha_1\geq\alpha_2\geq\ldots\geq\alpha_k>0$, 
$\alpha_i=0$ for $i>k$, and $\sum_{i=1}^k\alpha_i=n$. Usually, we will write $\alpha=(\alpha_1,\ldots,\alpha_k)$. 
The \textit{Young diagram} of $\alpha$ is an array of $n$
boxes, all of the same size and aligned to each other, with $k$ left-justified lines, 
where the $i$-th line contains $\alpha_i$ boxes, for all $i=1,\ldots,k$. Let $\beta=(\beta_1,\beta_2,\ldots)$ 
be an improper partition of $n$, i.e. a sequence of integers $\beta_i\geq 0$ such that $\beta_\ell>0$ and
$\beta_i=0$ for $i>\ell$ and $\sum_{i=1}^\ell\beta_i=n$. As for the proper partitions, 
we will write $\beta=(\beta_1,\ldots,\beta_\ell)$. A \textit{semistandard Young tableau
of shape $\alpha$ and of content $\beta$} 
(or \textit{generalized Young tableau of shape $\alpha$ and of content $\beta$})
is an array of numbers which is obtained from the Young diagram of $\alpha$ by inserting into
$\beta_i$ boxes the number $i$, for all $i$, such that 
\begin{itemize}
\item the entries in the rows of the diagram are increasing, and
\item the entries in the columns of the diagram are strictly increasing.
\end{itemize}
Usually one leaves the boxes away after having filled the diagram with numbers. 
Here is an example of a semistandard Young tableau of shape $\alpha=(4,4,3,3)$ and of content $\beta=(3,3,2,2,3,1)$:
\begin{equation*}
  \begin{array}{cccc}
    1 & 1 & 1 & 2 \\
    2 & 2 & 3 & 4 \\
    3 & 4 & 5 &   \\
    5 & 5 & 6 &   \\
  \end{array}
\end{equation*}

An important class of semistandard Young tableaux is the class of \textit{standard Young tableaux}. 
By a standard Young tableau of shape $\alpha$, we understand a semi\-standard Young tableau of shape $\alpha$ and of the particular content 
$\beta=\epsilon:=(1,\ldots,1)$, the partition consisting of $n$ times the number $1$. 
Of course, the entries of a standard Young tableau are strictly increasing in both the rows and the columns. 

Semistandard Young tableaux, and in particular standard Young tableaux, appear in various branches of mathematics, see the survey article \cite{sagan} 
for an overview. To these branches belong group representations, combinatorics, invariant theory, 
symmetric functions, the theory of algorithms and quantum algebras. However, in many of the instances where Young tableaux appear, 
a certain set of Young tableaux, subject to some constraint, has to be counted. The first counting formula for Young tableaux is due to 
Percy MacMahon (see \cite {macmahon}, vol. 2, sec. 429, or \cite{stanley}, p. 400), who used the language of plane partitions rather than  
of Young tableaux. From the wide range of more recent research, 
let me mention Donald Knuth's paper \cite{knuth}, where a correspondence between semistandard Young tableaux (of any shape and of any content) 
and symmetric matrices of nonnegative integers (of any size) is established. Paper \cite{knuth} served as an inspiration for several counting formulas for 
semistandard Young tableaux. We owe the first of these formulas to Basil Gordon \cite{gordon}. 
He studied the set $A_{p,q}$ of all semistandard Young tableaux with at most $q$ columns and with entries in $\{1,\ldots,p\}$, and proved that
\begin{equation}\label{apq}
 |A_{p,q}|=\prod_{1\leq i\leq j\leq p}\frac{q+i+j-1}{i+j-1}\,.
\end{equation}

Several refinements of this formula have so far appeared, let us discuss two of them in more detail. 
In each of these, the set of semistandard Young tableaux to be counted was shrunk by imposing more constraints on the shape of the tableaux. 
Myriam de Sainte-Catherine and G\'erard Viennot \cite{desainte} counted the number of elements of the set $B_{p,2q}$ 
of all semistandard Young tableaux with at most $2q$ columns, 
with an even number of elements in each row, and with entries in $\{1,\ldots,p\}$. 
Seul Hee Choi and Dominique Gouyou-Beauchamps \cite{choi} counted the number of elements of the set $C_{p,2q,r}$ 
of all semistandard Young tableaux with at most $2q$ columns, with an even number of elements in each row, with at most $r$ rows, 
and with entries in $\{1,\ldots,p\}$. Both counting formulas are structurally similar to \eqref{apq}. In particular, both are products. 

A list of further mathematicians who studied Young tableaux and enumeration formulas for these would include:
Ira Gessel and G{\'e}rard Viennot (\cite{gesselviennot}) George Andrews (e.g. \cite{andrews}, where another proof of 
Gordon's formula \eqref{apq} is given), Robert Proctor (\cite{proctor}), Christian Krattenthaler (\cite{krattenthaler1}, \cite{krattenthaler2}), 
and Ilse Fischer (\cite{fischer}).

The aim of the present paper is to count semistandard Young tableaux under the most restrictive constraints both on shape and on content:
We count semistandard Young tableaux of given shape $\alpha$ and of given content $\beta$. 
\begin{dfn}
  Let $\alpha$ be a proper partition of $n$ and $\beta$ an improper partition of $n$. Then define $K_{\alpha,\beta}$ to be the
  number of semistandard Young tableaux of shape $\alpha$ and of content $\beta$. 
\end{dfn}
The numbers $K_{\alpha,\beta}$ are called {\it Kostka numbers} in the literature (\cite{stanley}, p. 311).

\section{A classical motivation}

Originally, Alfred Young introduced his tableaux for studying representations of the symmetric group $S_n$. 
We will briefly discuss some of the outlines of his theory, since it serves well as a motivation 
for the theorem we are going to prove. It is a well-known fact that the inequivalent irreducible 
$\Q[S_n]$-modules correspond to the proper partitions of $n$. As reference, see e.g. \cite{jameskerber}, 
Theorem 2.1.11, or \cite{fulton}, Paragraph 7.2, Proposition 1. We will denote the irreducible module corresponding 
to the proper partition $\alpha$ of $n$ by $S^\alpha$. This module is called the \textit{Specht module corresponding to} $\alpha$. 
Let $f^{\alpha}$ be the number of standard Young tableaux of shape $\alpha$. Then a classical theorem states that
the $\Q$-dimension of $S^\alpha$ equals $f^\alpha$. As reference, see e.g. \cite{jameskerber}, Corollary 3.1.13. 

Another instance where the number $f^\alpha$ appears is the following situation: Consider the regular $\Q[S_n]$-module, 
i.e. $\Q[S_n]$ itself, viewed as a $\Q[S_n]$-module. One can decompose this module into its irreducible components. 
Now the question arises, how often does $S^\alpha$ appear in the decomposition of $\Q[S_n]$? 
The ans\-wer is: $f^\alpha$ times. On the one hand, this follows from the fact that $\Q[S_n]$ is isomorphic to a product of matrix algebras over $\Q$
(also to be found in \cite{jameskerber}), 
hence the $\Q$-dimension of each irreducible representation equals the number of times it appears in the regular representation. 
On the other hand, this is a special case of Young's rule. The full statement of Young's rule applies to a more general situation:  
Let $H$ be a subgroup of $S_n$. Then 
$\Q[S_n/H]=\oplus_{\overline{\sigma}\in S_n/H}\Q\overline{\sigma}$ is a $\Q[S_n]$-module with 
action of $\Q[S_n]$ on $\Q[S_n/H]$ induced by the action of $S_n$ on the cosets $S_n/H$. Consider the special case where the subgroup 
$H$ equals $H_\beta:=S_{\beta_1}\times\ldots\times S_{\beta_\ell}$ of $S_n$. Here we identify $S_{\beta_1}$ with the
symmetric group on $\{1,\ldots,\beta_1\}$, $S_{\beta_2}$ with the symmetric group on 
$\{\beta_1+1,\ldots,\beta_1+\beta_2\}$, etc. 
\begin{thm}[Young's rule, third version \cite{jameskerber}]
  The Specht module $S^\alpha$ appears in the decomposition of $\Q[S_n/H_\beta]$ 
  into irreducible $\Q[S_n]$-modules precisely $K_{\alpha,\beta}$ times. 
\end{thm}

There are two explicit formulas for $f^\alpha$. The first formula is the determinant formula
\begin{equation*}
  f^\alpha=n!\det(\frac{1}{(\alpha_i-i+j)!})
\end{equation*}
(where we define $1/{r!}=0$ whenever $r<0$) which, according to \cite{sagan}, goes back to Frobenius and Young. 
Expanding the determinant and using the definition of the multinomial coefficient, we can write this formula as
\begin{equation}\label{detf}
  f^\alpha=\sum_{\sigma\in S_k}{\rm sgn}(\sigma)\binom{n}{\alpha_1-1+\sigma(1),\ldots,\alpha_k-k+\sigma(k)}\,.
\end{equation}
(Remember that $k$ is the number of nonzero components of $\alpha$; as usual, $\binom{n}{r_1,\ldots,r_k}=0$ whenever some $r_{i}=0$.) 
The second formula is the hook formula 
\begin{equation}\label{hook}
  f^\alpha=\frac{n!}{\prod_{(i,j)\in\alpha}h_{(i,j)}}
\end{equation}
by Frame, Robinson and Thrall, as to be found in the original article \cite{frame}, or in \cite{fulton}, or in
\cite{jameskerber}. Here we used the coordinates $(i,j)$ in oder to locate the boxes of the diagram of $\alpha$.
Analogously as for matrices, $i$ numbers the lines downwards and $j$ numbers the columns from left to right. For a fixed $(i,j)$, 
the number $h_{(i,j)}=|H_{(i,j)}|$ is the length of the \textit{hook}
\begin{equation*}
  H_{(i,j)}=\{(i,j)\}\cup\{(i,j^\prime);j^\prime>j\}\cup\{(i^\prime,j);i^\prime>i\},
\end{equation*}
i.e. a subset of the diagram of $\alpha$ with ``knee'' in $(i,j)$, with ``leg" downwards and with ``arm" to the right. 
For both formulas, various proofs exist, some of them based on each other. The proofs use various ideas;
there are inductive, combinatorial and probabilistic proofs, see \cite{sagan} for an overview. 

The aim of this paper is to give an explicit formula for $K_{\alpha,\beta)}$ and to prove it by completely elementary means. 
Our formula will be a generalization of the determinant formula in its form \eqref{detf}. 
It turns out that all one has to do is to replace the multinomial coefficient in \eqref{detf} by another symmetric function, 
$\mu$, which we will define in Section \ref{secmu}. The technique of our proof here will be different from 
the techniques used in the proofs of $f^\alpha$, as listed in \cite{sagan}. Our proof will also be inductive, 
but will be based on the use of a certain functional equation. In Section \ref{secfunc}, we will 
prove that this functional equation is satisfied by the function $f$. In Section \ref{secthm}, we will state the theorem,
and give its proof, by showing that our formula for $K_{\alpha,\beta)}$ defines a function that satisfies the same functional equation as $f$. 

\section{A generalization of the multinomial coefficient}\label{secmu}

Let $E$ be the set of all finitely supported sequences with integer values. Let $s:E\to\Z$ be the 
function that sums the components of an element of $E$; we will need the element $0:=(0,0,\ldots)$ of $E$; we will need, 
for all $k\in\N$, the element $(k):=(1,2,\ldots,k,0,\ldots)$ of $E$ and we will need, for all $k\in\N$ and for all $\sigma\in S_k$, 
the element $(\sigma(k)):=(\sigma(1),\sigma(2),\ldots,\sigma(k),0,\ldots)$ of $E$. 
We will also need sums and differences of elements of $E$. In other words, $E$ will have to carry a $\Z$-module structure.
For our purposes, the right way to define addition and scalar multiplication is the most natural way, i.e. component\-wise.
Further, for every sequence $\rho=(\rho_1,\rho_2,\ldots)$ of integers, we define $E_\rho$ 
to be the set of all elements $\delta$ of $E$ such that $s(\delta)=0$, 
or $s(\delta)=\rho_1$, or $s(\delta)=\rho_1+\rho_2$, etc. Finally, let $E_+$ be the set of sequences in $E$
with nonnegative values. 
\begin{dfn}\label{mu}
  Define a map $\mu_\rho:E_\rho\to\N_0$ by the three properties
  \begin{itemize}
    \item $\mu_\rho(0)=1$,
    \item $\mu_\rho(\delta)=0$ if some $\delta_i<0$,
    \item $\mu_\rho(\delta)=\sum_{\gamma\in E_+,s(\gamma)=\rho_\ell}\mu_\rho(\delta-\gamma)$ 
    if $s(\delta)=\rho_1+\ldots+\rho_\ell$. 
  \end{itemize}
  Extend the domain of definition of $\mu_\rho$ by setting $\mu_\rho(\delta):=0$ for all $\delta\notin E_\rho$. 
\end{dfn}

Let us sketch some of the values of  $\mu_\rho$ in example $\rho=(3,2,3,2,\ldots)$. We compute some 
of the values of $\mu_\rho(\delta)$ for sequences of the form $\delta=(\delta_1,\delta_2,0,\ldots)$. 
In the following matrix, we insert $\mu_\rho(\delta)$ at the position $(\delta_1,\delta_2)$. 
\begin{equation*}
  \begin{array}{ccccccccc}
    1 & 0 & 0  & 1  & 0  & 1  & 0  & 0   & 1   \\
    0 & 0 & 1  & 0  & 2  & 0  & 0  & 3   & 0   \\
    0 & 1 & 0  & 3  & 0  & 0  & 6  & 0   & 10  \\
    1 & 0 & 3  & 0  & 0  & 9  & 0  & 18  & 0   \\
    0 & 2 & 0  & 0  & 10 & 0  & 25 & 0   & 0   \\
    1 & 0 & 0  & 9  & 0  & 28 & 0  & 0   & 81  \\
    0 & 0 & 6  & 0  & 25 & 0  & 0  & 96  & 0   \\
    0 & 3 & 0  & 18 & 0  & 0  & 96 & 0   & 273 \\
    1 & 0 & 10 & 0  & 0  & 81 & 0  & 273 & 0   \\ 
  \end{array}
\end{equation*}
The matrix should be thought of as having infinite size. Towards the left and upwards, all entries of the matrix are 
zero, since $\mu_\rho(\delta)=0$ if some $\delta_i<0$. 

The map $\mu_\rho$ is indeed a generalization of the multinomial coefficient. More precisely, 
if $\delta=(\delta_1,\ldots,\delta_k)$ and $s(\delta)=n$, then 
\begin{equation*}
  \binom{n}{\delta_1,\ldots,\delta_k}=\mu_\epsilon(\delta)
\end{equation*}
for $\epsilon=(1,\ldots,1)$ as before. This is a consequence of the following lemma. 
\begin{lmm}
  Let $X_1,\ldots,X_k$ be indeterminates over $\Z$, and let $\delta=(\delta_1,\ldots,\delta_k,0,\ldots)$ be a 
  sequence in $E_\rho$ such that $s(\delta)=\rho_1+\ldots+\rho_\ell$. Then the number $\mu_\rho(\delta)$ 
  equals the coefficient of the monomial $X_1^{\delta_1}\ldots X_k^{\delta_k}$ when expanding the product
  \begin{equation*}
    (\sum_{e_{1,1}+\ldots+e_{1,k}=\rho_1}X_1^{e_{1,1}}\ldots X_k^{e_{1,k}})\ldots(\sum_{e_{\ell,1}+\ldots+e_{\ell,k}
    =\rho_\ell}X_1^{e_{\ell,1}}\ldots X_k^{e_{\ell,k}})\,,
  \end{equation*}
  where the exponents $e_{i,j}$ run through $\N_{0}$. 
\end{lmm} 

The proof of this lemma is straightforward induction over $\ell$.
Note that the lemma implies in particular that $\mu_\rho$ is a symmetric function in the entries of $\delta$. 

We defined the function $\mu_\rho$ recursively. In order to establish a connection between $K_{\alpha,\beta}$ 
and $\mu_\beta$, we also need some recursion for $K_{\alpha,\beta}$. This is done in the following section. 

\section{A functional equation for $f$}\label{secfunc}

We would like to determine the number $K_{\alpha,\beta}$, which counts how many ways there are of 
filling the Young diagram of $\alpha$ with numbers $1,\ldots,\ell$ such that the result is a semistandard Young tableau of shape 
$\alpha$ and of content $\beta$. Let us try to fill the Young diagram of $\alpha$ little by little, starting with large
numbers and then going downwards. In this paragraph we will make the first and decisive step: We will 
insert the highest number, i.e. $\ell$, precisely $\beta_{\ell}$ times into the Young diagram of $\alpha$. 
From this step, we will derive the desired functional equation for $f$. 

The obvious question is: Which possibilities are there of inserting the number $\ell$ precisely $\beta_{\ell}$ times into the 
Young diagram of $\alpha$ such that the boxes remaining unfilled can be filled with 
$\beta_{1}$ times $1$, $\beta_{2}$ times $2$, and so on, all the way up to $\beta_{\ell-1}$ times $\ell-1$,
yielding a semistandard Young tableau of shape $\alpha$ and of content $\beta$?

In order to answer this question, let us take a look at a potential result of the filling process. More concretely, let us take a look at a fixed
semistandard Young tableau of shape $\alpha$ and of content $\beta$. The entries in the rows of the semistandard Young tableau
are increasing. This imposes the following necessary condition on the position of $\ell$: No lesser number than $\ell$ may be inserted
to the right of an $\ell$ in any row of the semistandard Young tableau. In other words, the numbers $\ell$ appearing in a fixed row 
of the semistandard Young diagram are aligned at the right-hand end of the row. For all $i\in\{1,\ldots,k\}$, let $\gamma_{i}$
denote the number of occurrences of $\ell$ in the $i$-th row of the semistandard Young tableau. 
Put $\gamma=(\gamma_1,\gamma_2,\ldots)$. Then $\gamma$ is clearly an element of $E_{+}$, and $s(\gamma)=\beta_\ell$. 

So far we have used the constraint on the rows of the semistandard Young tableau. The constraint on the columns 
of the semistandard Young tableau implies a constraint on $\gamma$, which we formulate as a lemma. 
\begin{lmm}\label{lmmgamma}
  For all $i\in\{1,\ldots,k\}$, we have 
  \begin{equation}\label{gamma}
    \alpha_i-\gamma_i\geq\alpha_{i+1}\,.
  \end{equation}
\end{lmm}

\begin{proof}
  Let us take a closer look at the lines number $i$ and number $i+1$ of our fixed semistandard Young tableau
  of shape $\alpha$ and of content $\beta$. Since we are interested only in the question where $\ell$ is inserted, 
  we replace every number different from $\ell$ by an asterisk. The result looks as follows:
  \begin{equation*}
    \begin{array}{ccccccc}
      * & * & * & * & *      & \ell & \ell \\ 
      * & * & * & * & \ell &        &      
    \end{array}
  \end{equation*}
  Since the entries of a semistandard Young tableau are strictly increasing, a picture like 
  \begin{equation*}
    \begin{array}{ccccccc}
      * & * & * & * & \ell & \ell & \ell \\ 
      * & * & * & * & *      &        &      
    \end{array}
  \end{equation*}
  or 
  \begin{equation*}
    \begin{array}{ccccccc}
      * & * & * & * & \ell & \ell & \ell \\ 
      * & * & * & * & \ell &        &      
    \end{array}
  \end{equation*}
  does not occur. The abstract meaning of these three pictures is that the number of asterisks in the $i$-th row may not exceed 
  the number of entries in the $(i+1)$-st row. Equation \eqref{gamma} translates this fact to a formula. 
\end{proof}

It is important to note that \eqref{gamma} is equivalent to $\alpha-\gamma$ being a proper partition of $n-\beta_\ell$. 

\begin{pro}
  Let $\alpha$ and $\beta$ be defined as above. Let $\beta^{\prime}:=(\beta_1,\ldots,\beta_{\ell-1},0,\ldots)$. 
  Then the function $f$ is subject to the following functional equation:
  \begin{equation}\label{functionalf}
    K_{\alpha,\beta}=\sum_{\gamma\in E_+,s(\gamma)=\beta_\ell\,,\text{\eqref{gamma} holds}}
    K_{\alpha-\gamma,\beta^\prime}\,.
  \end{equation}
\end{pro}

\begin{proof}
  Lemma \eqref{lmmgamma} tells us that condition \eqref{gamma} is necessary in the following sense:
  For every semistandard Young tableau of shape $\alpha$ and of content $\beta$, there exists a $\gamma\in E_{+}$ satisfying \eqref{gamma}
  such that the given semistandard Young tableau of shape $\alpha$ and of content $\beta$ is obtained from
  a semistandard Young tableau of shape $\alpha-\gamma$ and of content $\beta^{\prime}$ by the following process:
  Take the semistandard Young tableau of shape $\alpha-\gamma$ and of content $\beta^{\prime}$ and append 
  $\gamma_{i}$ times the number $\ell$ to the $i$-th row, for all $i\in\{1,\ldots,k\}$. 
  It is clear that for a given semistandard Young tableau of shape $\alpha$ and of content $\beta$, the 
  semistandard Young tableau of shape $\alpha-\gamma$ and of content $\beta^{\prime}$ which under this process leads to the 
  given semistandard Young tableau of shape $\alpha$ and of content $\beta$, is unique (since, conversely, the 
  semistandard Young tableau of shape $\alpha-\gamma$ and of content $\beta^{\prime}$ is obtained by erasing 
  $\gamma_{i}$ times the number $\ell$ in every line).
  
  In an analogous sense, condition \eqref{gamma} is also sufficient. Let us express sufficiency as follows:
  Suppose to be given, along with $\alpha$ and $\beta$, an element $\gamma$ of $E_{+}$ satisfying \eqref{gamma}. 
  Then for every semistandard Young tableau of shape $\alpha-\gamma$ and of content $\beta^{\prime}$, 
  there exists a unique semistandard Young tableau of shape $\alpha$ and of content $\beta$ which is obtained from the given 
  semistandard Young tableau of shape $\alpha-\gamma$ and of content $\beta^{\prime}$ by appending
  $\gamma_{i}$ times the number $\ell$ to the $i$-th row, for all $i\in\{1,\ldots,k\}$. 
  
  Hence there is a bijection between the set of semistandard Young tableaux of shape $\alpha$ and of content $\beta$,
  and the set of semistandard Young tableaux of shape $\alpha-\gamma$ and of content $\beta^{\prime}$ such that 
  $\gamma$ satisfies \eqref{gamma}. Counting both sets gives functional equation \eqref{functionalf}. 
\end{proof}

\section{The counting formula}\label{secthm}

Before stating the theorem, let us rule out one potential obstruction. 
A given improper partition $\beta=(\beta_1,\ldots,\beta_\ell)$ of $n$ may contain some components $\beta_i=0$, 
for $i<\ell$. At first sight, this seems to cause difficulty when determining $K_{\alpha,\beta}$. Yet, it does not do so, 
since we can get rid of the ``gaps" in $\beta$ by the following process: 
Given $\beta$, we remove from $\beta$ all $\beta_i=0$, where $i<\ell$, push the remaining components of $\beta$ to the left 
and call the result $\rho$. This is a sequence with the same nonzero components as $\beta$, and the nonzero components of $\rho$ 
appear in the same order as the nonzero components of $\beta$. For the time being, let $B$ denote the set of 
subscripts of the nonzero elements of $\beta$, let $R$ denote the set of subscripts of the nonzero elements of $\rho$, 
and let $t$ denote the unique strictly monotonous bijection $t:R\to B$. Then the entries of a semistandard Young tableau 
of shape $\alpha$ and of content $\beta$ are clearly $\beta_{t(1)}$ times the number $t(1)$, $\beta_{t(2)}$ times the number $t(2)$, and so on, 
all the way up to $\beta_{t(|R|)}$ times the number $t(|R|)$. Given a semistandard Young tableau of shape $\alpha$ and of content $\beta$, 
we can replace each $t(i)$ by $i$. The result will be a semistandard Young tableau of shape $\alpha$ and of content $\rho$. 
(The fact that the resulting array of numbers satisfies the monotony conditions for a semistandard Young tableau follows from the 
strict monotony of $t$.) Conversely, given a semistandard Young tableau of shape $\alpha$ and of content $\rho$, 
we can replace each $i$ by $t(i)$ and will thus obtain a semistandard Young tableau of shape $\alpha$ and of content $\beta$. 
(Again, we use the strict monotony of $t$.) Thus $K_{\alpha,\beta}=K_{\alpha,\rho}$. 
Therefore, for the rest of the paper, we can assume w.l.o.g. that the improper partition 
$\beta=(\beta_1,\ldots,\beta_\ell)$ contains no components $\beta_i=0$, for $i\leq\ell$. 

\begin{thm}\label{it}
  The number $K_{\alpha,\beta}$ of semistandard Young tableaux of shape $\alpha$ and of content $\beta$ is given by 
  \begin{equation}\label{thmform}
    K_{\alpha,\beta}=\sum_{\sigma\in S_k}{\rm{sgn}}(\sigma)\mu_\beta\bigl(\alpha-(k)+(\sigma(k))\bigr)\,.
  \end{equation}
\end{thm}

\begin{proof}
  Let us denote the function on the right-hand side of \eqref{thmform} (a function taking arguments $\alpha$ and $\beta$) by $g_{\alpha,\beta}$. 
  Obviously $g_{\alpha,\beta}$ is defined for all $(\alpha,\beta)\in E\times E$, 
  unlike $K_{\alpha,\beta}$, which is defined only for those $(\alpha,\beta)$ in $E\times E$
  for which $\alpha$ is a proper partition of $n$ and $\beta$ is an improper partition of $n$, 
  for some $n\leq 1$. We have to prove that $K_{\alpha,\beta}=g_{\alpha,\beta}$ for all $(\alpha,\beta)$ in
  the domain of definition of $f$. 
  We will prove the theorem by first showing that $g$ also satisfies functional equation \eqref{functionalf}
  and then showing that $f$ and $g$ satisfy the same boundary condition, in a sense that will be explained more precisely later. 

  So let us start with the functional equation for $g$. From the definition of function $\mu_{\beta}$ it follows immediately that 
  \begin{equation*}
    g_{\alpha,\beta}=\sum_{\gamma\in E_+,s(\gamma)=\beta_\ell}g_{\alpha-\gamma,\beta}\,.
  \end{equation*}
  When computing $g_{\alpha-\gamma,\beta}$, 
  we sum over certain $\mu_\beta(a)$, where $a\in E_{\beta}$ and $s(a)=\beta_1+\ldots+\beta_{\ell-1}$. 
  From the definition of function $\mu_{\beta}$ follows that for computing of $\mu_\beta(a)$, where $s(a)=\beta_1+\ldots+\beta_{\ell-1}$, 
  only the first $\ell-1$ terms of $\beta$ are relevant. As before, let us write $\beta^{\prime}=(\beta_1,\ldots,\beta_{\ell-1},0,\ldots)$. 
  Then clearly $\mu_\beta(a)=\mu_{\beta^\prime}(a)$ for all $a\in E_{\beta}$ with $s(a)=\beta_1+\ldots+\beta_{\ell-1}$. 
  This implies that 
  \begin{equation}\label{functionalg}
    g_{\alpha,\beta}=\sum_{\gamma\in E_+,s(\gamma)=\beta_\ell}g_{\alpha-\gamma,\beta^\prime}\,.
  \end{equation}
  
  Equation \eqref{functionalg} is already very similar to equation \eqref{functionalf}. 
  The difference is that in \eqref{functionalg}, sum is taken over a larger set than in \eqref{functionalf}. 
  Let us denote by $X$ the difference between the indexing sets of the sums in \eqref{functionalg} and in \eqref{functionalf}. 
  Thus $X$ is the set of all $\gamma\in E_+$ such that $s(\gamma)=\beta_\ell$, but \eqref{gamma} does not hold. 
  We will now prove the following assertion: 
  \begin{equation}\label{not}
    \sum_{\gamma\in X}g_{\alpha-\gamma,\beta^\prime}=0\,.
  \end{equation}
  
  Let us define a map $\xi:X\to X$ in the following way: Given $\gamma\in X$, look for the smallest $i$ 
  such that $\alpha_i-\gamma_i<\alpha_{i+1}$. Then set 
  \begin{equation*}
    \xi(\gamma)_j=
    \begin{cases}
    \gamma_j & \mbox{for} \quad j\neq i,i+1\,,\\ 
   \alpha_i-\alpha_{i+1}+\gamma_{i+1}+1 & \mbox{for} \quad j=i\,,\\ 
   \alpha_{i+1}-\alpha_i+\gamma_i-1 & \mbox{for} \quad j=i+1\,.
  \end{cases}
  \end{equation*}
  It is clear that $\xi(\gamma)$ lies in $E_{+}$ and that $s(\xi(\gamma))=\beta_{\ell}$. 
  Since $\alpha_i-\xi(\gamma)_i=\alpha_i-(\alpha_i-\alpha_{i+1}+\gamma_{i+1}+1)<\alpha_{i+1}$, 
  sequence $\xi(\gamma)$ does not satisfy \eqref{gamma}. Thus $\xi:X\to X$ is indeed well defined. 

  Now $\gamma$ and $\xi(\gamma)$ differ only in the $i$-th and in the $(i+1)$-st component. 
  For the forthcoming discussion, let us fix the transposition $\tau=(i\,\,i+1)$. 
  Let us compare the $i$-th and the $(i+1)$-st component of $\alpha-\gamma-(k)+(\sigma(k))$ 
  and $\alpha-\xi(\gamma)-(k)+(\sigma\circ\tau(k))$. On the one hand, we have 
  \begin{equation}\label{1}
    \begin{split}
      &\bigl(\alpha-\gamma-(k)+(\sigma(k))\bigr)_i=\alpha_i-\gamma_i-i+\sigma(i)\\
      = & \alpha_{i+1}-\xi(\gamma)_{i+1}-(i+1)+\sigma(i)=\bigl(\alpha-\xi(\gamma)-(k)+(\sigma\circ\tau(k))\bigr)_{i+1}\,,\\
    \end{split}
  \end{equation}
  and on the other hand, we have
  \begin{equation}\label{2}
    \begin{split}
      &\bigl(\alpha-\gamma-(k)+(\sigma(k))\bigr)_{i+1}=\alpha_{i+1}-\gamma_{i+1}-(i+1)+\sigma(i+1)\\
      = & \alpha_i-\xi(\gamma)_i-i+\sigma(i+1)=\bigl(\alpha-\xi(\gamma)-(k)+(\sigma\circ\tau(k))\bigr)_i\,.
    \end{split}
  \end{equation}
  Further, for all $j\neq i,i+1$, we clearly have
  \begin{equation}\label{3}
    \bigl(\alpha-\gamma-(k)+(\sigma(k))\bigr)_{j}=\bigl(\alpha-\xi(\gamma)-(k)+(\sigma\circ\tau(k))\bigr)_{j}.
  \end{equation}
  Equations \eqref{1}, \eqref{2}, \eqref{3}, together with the fact that $\mu_\beta$ is a symmetric function, imply that 
  \begin{equation}\label{symm}
    \mu_{\beta^{\prime}}\bigl(\alpha-\gamma-(k)+(\sigma(k))\bigr)
    =\mu_{\beta^{\prime}}\bigl(\alpha-\xi(\gamma)-(k)+(\sigma\circ\tau(k))\bigr)\,.
  \end{equation}
  The map $S_{k}\to S_{k}:\sigma\mapsto\sigma\circ\tau$ is a bijection. Clearly, ${\rm{sgn}}(\sigma)=-{\rm{sgn}}(\sigma\circ\tau)$. 
  Therefore \eqref{symm}, along with the definition of function $g$, yields 
  $g_{\alpha-\gamma,\beta^\prime}=-g_{\alpha-\xi(\gamma),\beta^\prime}$ for the case that 
  the partition $\alpha-\gamma$ consists of $k$ nonzero parts. If the partition $\alpha-\gamma$ consists of less than $k$ nonzero parts, 
  the above discussion translates literally to the situation where every $k$ is replaced by the number of nonzero parts in the partition $\alpha-\gamma$. 
  Thus $g_{\alpha-\gamma,\beta^\prime}=-g_{\alpha-\xi(\gamma),\beta^\prime}$ for all $\gamma\in X$. 
  Since $\xi:X\to X$ is a bijection, it follows that 
  $\sum_{\gamma\in X}g_{\alpha-\gamma,\beta^\prime}=-\sum_{\gamma\in X}g_{\alpha-\gamma,\beta^\prime}$, 
  hence $\sum_{\gamma\in X}g_{\alpha-\gamma,\beta^\prime}=0$, as claimed. 

  Therefore, function $g$ also satisfies functional equation \eqref{functionalf}. Applying the functional equation several times, for both $f$ and $g$,
  we finally arrive at a point where the improper partition $\beta$ in the second argument has only got one nonzero component, i.e., 
  $\beta=(\beta_1,0,\ldots)$. Hence it suffices to show that for this particular $\beta$, we have $K_{\alpha,\beta}=g_{\alpha,\beta}$, 
  for all proper partitions $\alpha$. This is the boundary condition for $f$ and $g$, announced already at the beginning of the proof.
  Therefore, for the rest of the proof we make the following assumptions: $n$ is arbitrary, $\alpha=(\alpha_{1},\ldots,\alpha_{k})$ 
  is a proper partition of $n$ and $\beta=(n,0,\ldots)$. 
  
  If we try to construct semistandard Young tableaux of shape $\alpha$ and of content $\beta$, we have to fill the 
  Young diagram of $\alpha$ with $n$ times the number $1$ such that, in particular, the entries 
  in every column are strictly increasing. This is only possible if the Young diagram of $\alpha$ has just one row, i.e., if $\alpha=\beta$. 
  In this case, there is a unique semistandard Young tableau of shape $\alpha$ and of content $\beta$. 
  Thus $K_{\beta,\beta}=1$, and all $K_{\alpha,\beta}=0$. where $\alpha\neq\beta$. 
  In other words, $K_{\beta,\beta}=1$ if $k=1$ and $K_{\alpha,\beta}=0$ if $k\geq2$. We have to show that this is also true for $g$. 

  Let us treat the two cases $k=1$ and $k\geq2$ separately. 
  For $k=1$, the only possibility for $\alpha$ to be a proper partition of $n$ is the case $\alpha=(n,0,\ldots)$, i.e., 
  $\alpha=\beta$. In this case, the only summand occurring in the sum defining $g_{\alpha,\beta}$ 
  is the summand for $\sigma={\rm{id}}$, and this summand yields $g_{\alpha,\beta}=\mu_\beta(\beta)=1$
  by definition of $\mu_\beta$. For $k\geq2$, we have to show that $g_{\alpha,\beta}=0$. 
  In order to determine $g_{\alpha,\beta}$, we have to determine $\mu_{\beta}\bigl(\alpha-(k)+(\sigma(k))\bigr)$ for
  all $\sigma\in S_{k}$. 
  For a given $\sigma\in S_{k}$, we distinguish between the following two cases: Either there is some $i$ such that $\alpha_i-i+\sigma(i)<0$,  
  in this case $\mu_{\beta}\bigl(\alpha-(k)+(\sigma(k))\bigr)=0$, or all $\alpha_i-i+\sigma(i)\geq0$, then
$\mu_{\beta}\bigl(\alpha-(k)+(\sigma(k))\bigr)=1$. 
  Let $Y$ be the set of those $\sigma\in S_{k}$ for which all $\alpha_i-i+\sigma(i)\geq0$. 
  Then for all $\sigma\in Y$, the summand in $g_{\alpha,\beta}$ corresponding to $\sigma$ equals ${\rm{sgn}}(\sigma)$. 
  We thus have to show that $\sum_{\sigma\in Y}{\rm{sgn}}(\sigma)=0$. 
  Let us do this in a way analogous to what we have done before, namely by making use of an appropriate transposition. 
  Here it is going to be the fixed transposition $\tau=(1\,\,2)$. We proceed as follows: 
  Our first observation is that from $\alpha$ being a proper partition of $n$ and $k\geq2$, we get in particular that $\alpha_1\geq\alpha_2\geq1$. 
  From that we deduce that $\alpha_1-1+\sigma(1)\geq1$ and $\alpha_2-2+\sigma(2)\geq0$ for all $\sigma\in S_{k}$. 
  Further, we also have $\alpha_1-1+\sigma(2)\geq1$ and $\alpha_2-2+\sigma(1)\geq0$ for all $\sigma\in S_{k}$. 
  In particular the first two components of $\alpha-(k)+(\sigma(k))$ and $\alpha-(k)+(\sigma\circ\tau(k))$ are nonnegative. 
  And clearly, $\alpha-(k)+(\sigma(k))$ and $\alpha-(k)+(\sigma\circ\tau(k))$ differ only in the first two components. 
  Thus, in particular, $\alpha-(k)+(\sigma(k))$ and $\alpha-(k)+(\sigma\circ\tau(k))$ have the same negative components. 
  Now $\mu_\beta\bigl(\alpha-(k)+(\sigma(k))\bigr)$ and $\mu_\beta\bigl(\alpha-(k)+(\sigma\circ\tau(k))\bigr)$ can both
  only take the values $0$ or $1$, depending on whether they have some negative component or not. Therefore, 
  \begin{equation}\label{invariance}
    \mu_\beta\bigl(\alpha-(k)+(\sigma(k))\bigr)=\mu_\beta\bigl(\alpha-(k)+(\sigma\circ\tau(k))\bigr)\,. 
  \end{equation}
  Since the map $S_{k}\to S_{k}:\sigma\mapsto\sigma\circ\tau$ is a bijection, 
  and ${\rm{sgn}}(\sigma)=-{\rm{sgn}}(\sigma\circ\tau)$, equation \eqref{invariance} yields 
  $\sum_{\gamma\in S_k}{\rm{sgn}}(\sigma)\mu_\beta\bigl(\alpha-(k)+(\sigma(k))\bigr)
  =-\sum_{\gamma\in S_k}{\rm{sgn}}(\sigma)\mu_\beta\bigl(\alpha-(k)+(\sigma(k))\bigr)$, 
  hence $g_{\alpha,\beta}=0$, as claimed. 
\end{proof}

\section{Comments on the literature}\label{literature}

After having written a first draft of this paper, I was told by Mark Shimozono and Christian Krattenthaler that Theorem \ref{it}
follows from classical results. Here is Christian Krattenthaler's sketch of a proof. 

First observe that the number $\mu_\rho(\delta)$ equals the number of $\N_{0}$-matrices $A$ satisfying ${\rm row}(A)=\rho$ and 
${\rm col}(A)=\delta$, where 
\begin{equation*}
  {\rm row}(A)=(\sum_{j}a_{1,j},\sum_{j}a_{2,j},\ldots)\,, \text{ and }
  {\rm col}(A)=(\sum_{i}a_{i,1},\sum_{i}a_{i,2},\ldots)\,.
\end{equation*}
The number of such matrices is denoted $N_{\rho,\delta}$ (see e.g. \cite{stanley}, Proposition 7.5.1). Let us consider $N$ as 
a function on the set of pairs $(\rho,\delta)$ of elements of $E$ such that $s(\rho)=s(\delta)$. This constraint on the arguments of $N$ 
reflects the fact that for every matrix $A$, we have $s({\rm row}(A))=s({\rm col}(A))$. One shows 
that $N_{\rho,\delta}=\mu_\rho(\delta)$ for all $(\rho,\delta)$ in the domain of definition of $N$ by verifying that for a fixed $\rho$, 
the function $E_{\rho}\to\N_{0}:\delta\mapsto N_{\rho,\delta}$ satisfies the three items from Definition \ref{mu}. 
The first two are trivial; the third reads 
\begin{equation}\label{N}
  N_{\rho,\delta}=\sum_{\gamma\in E_+,s(\gamma)=\rho_\ell}N_{\rho^\prime,\delta-\gamma}\,,\text{ if }
  s(\delta)=\rho_1+\ldots+\rho_\ell\,,
\end{equation}
where $\rho^\prime=(\rho_{1},\ldots,\rho_{\ell-1})$. 
Note that an equation analogous to \eqref{N} also holds for $\mu_\rho(\delta)$ and $\mu_{\rho^\prime}(\delta-\gamma)$ instead of
$N_{\rho,\delta}$ and $N_{\rho^\prime,\delta-\gamma}$, as has been remarked in the proof of Theorem \ref{it}. 
Now for proving \eqref{N}, take a matrix $A$ such that ${\rm row}(A)=\rho$ and ${\rm col}(A)=\delta$ and delete its $\ell$-th row. 
The resulting matrix, say $B$, will satisfy ${\rm row}(B)=\rho^\prime$ and ${\rm col}(B)=\delta-\gamma$, for a uniquely determined 
$\gamma\in E_+$ such that $s(\gamma)=\rho_\ell$. Conversely, to every $B$ satisfying ${\rm row}(B)=\rho^\prime$ and 
${\rm col}(B)=\delta-\gamma$, for some $\gamma\in E_+$ such that $s(\gamma)=\rho_\ell$, can be added an $\ell$-th row such that 
the resulting matrix, say $A$, satisfies ${\rm row}(A)=\rho$ and ${\rm col}(A)=\delta$, and the $\ell$-th row thus added is unique. 
This proves \eqref{N}, hence $N_{\rho,\delta}=\mu_\rho(\delta)$.

Next, let us cite some facts from \cite{stanley}. We need the {\it monomial symmetric functions} $m_{\beta}$ (\cite{stanley}, 7.3), 
the {\it complete homogeneous symmetric functions} $h_{\lambda}$ (\cite{stanley}, 7.5), and the {\it Schur functions} $s_{\alpha}$
(\cite{stanley}, 7.10), which are defined for partitions $\beta$, $\lambda$, $\alpha$ of the integer $n$. We do not need the definitions 
of any of these functions, let us just remark that they are power series in indeterminates $X_{1},X_{2},\ldots$, and that the definition 
of $s_{\alpha}$ involves the set of all semistandard Young tableaux of shape $\alpha$. However, we do need the fact that the various 
$m_{\beta}$ form a basis of the $\Q$-vector space of {\it symmetric functions} over $\Q$, and we need the following identities:
\begin{equation}\label{h1}
  h_{\lambda}=h_{\lambda_{1}}h_{\lambda_{2}}\ldots\text{ if }\lambda=(\lambda_{1},\lambda_{2},\ldots)
\end{equation}
(this is the definition of $h_{\lambda}$ in 7.5),
\begin{equation}\label{h2}
  h_{\lambda}=\sum_{\beta}N_{\lambda,\beta}m_{\beta}\,,
\end{equation}\label{s}
where the sum ranges over all proper partitions of $n$ (this is Proposition 7.5.1),
\begin{equation}\label{m}
  s_{\alpha}=\sum_{\beta}K_{\alpha,\beta}m_{\beta}\,,
\end{equation}
where the sum again ranges over all proper partitions of $n$ (this follows from Definition 7.10.1), and
\begin{equation}\label{det}
  s_{\alpha}=\det(h_{\alpha_{i}-i+j})_{i,j=1}^{k}\,,
\end{equation}
where $\alpha=(\alpha_1,\ldots,\alpha_k)$ (this is Theorem 7.16.1). 

Now expanding the determinant in \eqref{det}, using \eqref{h1} for
$\lambda=\alpha-(k)+(\sigma(k))$, expressing all $h_{\alpha-(k)+(\sigma(k))}$ as linear combinations of the basis elements $\beta$ by 
\eqref{h2}, and so doing also on the left hand side of \eqref{det} by \eqref{m}, we obtain the formula from Theorem \ref{it}. 

Mark Shimozono has remarked that the formula from Theorem \ref{it} is the special case for general linear groups of a formula for the
dimension of the $\beta$ weight-space in the irreducible highest weight module of highest weight $\alpha$. Such a formula holds for any
weight space in any finite-dimensional highest-weight module over a simple Lie group and is due to Bertram Kostant, see \cite{fultonharris}, 
p. 419--424, or the original article \cite{kostant}. The formula is derived from Weyl's character formula (\cite{fultonharris}, p. 399--414). 
Finally, it is well known that the above dimension is equal to the number of semistandard Young tableaux of shape $\alpha$ and of content $\beta$, 
see \cite{fulton}, p. 121. 

These two proofs of Theorem \ref{it} are apparantly much shorter than the one given in the previous section. However, there is a large amount
of work behind the small number of theorems from which, as I have outlined in this section, Theorem \ref{it} can be derived in a quick and easy 
way. To my mind, the very value of the proof given in the previous section lies in its avoiding of any theoretical apparatus whatsoever.

\section{Concluding remarks}

In the first section, we cited Gordon's counting formula \eqref{apq} concerning the size of the set $A_{p,q}$ of all
semistandard Young tableaux with at most $q$ columns, and with entries in $\{1,\ldots,p\}$. So let us fix natural numbers
$p$ and $q$, and let $Z$ be the set of pairs $(\alpha,\beta)$, where $\alpha=(\alpha_{1},\ldots,\alpha_{k})$ is a proper partition of some
number $n$ such that all $\alpha_{i}\leq q$, and $\beta=(\beta_{1},\ldots,\beta_{\ell})$ is an improper partition of
the same $n$ such that $\ell\leq q$. Then clearly
\begin{equation*}
  \prod_{1\leq i\leq j\leq p}\frac{q+i+j-1}{i+j-1}=\sum_{(\alpha,\beta)\in Z}K_{\alpha,\beta}\,.
\end{equation*}
Together with formula \eqref{thmform} for the summands $K_{\alpha,\beta}$ on the right hand side, 
this is a nontrivial equality. It is fair to ask for a conceptual explanation of this equality. (In the same way, we 
get nontrivial equalities when we replace $A_{p,q}$ by $B_{p,2q}$, or by $C_{p,2q,r}$. Here, the question for an
explanation of the respective equalities remains open.)

Further questions concern the integration of Theorem \ref{it} into Mackey theory, and, of course, 
a generalization of the hook formula for Young tableaux of given shape and of given content. 

\section{Acknowledgments}

I did most of the work on this paper when writing my Ph.D. Thesis at the University of Innsbruck with Kurt Girstmair, whom I
want to thank for his advice. I also want to thank Gordon James and John McKay for their encouragement, Thomas Zink for 
very useful comments on this text, and Mark Shimozono and Christian Krattenthaler for communicating what I am trying to reproduce in 
Section \ref{literature}.

Many thanks to my brother Thomas Lederer for helping to improve my English and to Sarah L\"ow for technical support.

\bibliography{littableaux}

\providecommand{\bysame}{\leavevmode\hbox to3em{\hrulefill}\thinspace}
\providecommand{\MR}{\relax\ifhmode\unskip\space\fi MR }
\providecommand{\MRhref}[2]{%
  \href{http://www.ams.org/mathscinet-getitem?mr=#1}{#2}
}
\providecommand{\href}[2]{#2}
\begin{thebibliography}{10}

\bibitem{andrews}
George~E. Andrews, \emph{Plane partitions. {II}. {T}he equivalence of the
  {B}ender-{K}nuth and {M}ac{M}ahon conjectures}, Pacific J. Math. \textbf{72}
  (1977), no.~2, 283--291. \MR{MR472108 (80d:05006b)}

\bibitem{choi}
Seul~Hee Choi and Dominique Gouyou-Beauchamps, \emph{Enumeration of generalized
  {Y}oung tableaux with bounded height}, Theoret. Comput. Sci. \textbf{117}
  (1993), no.~1-2, 137--151. \MR{MR1235174 (94h:05093)}

\bibitem{desainte}
Myriam de~Sainte-Catherine and G{\'e}rard Viennot, \emph{Enumeration of certain
  {Y}oung tableaux with bounded height}, Combinatoire \'enum\'erative
  (Montreal, Que., 1985/Quebec, Que., 1985), Lecture Notes in Math., vol. 1234,
  Springer, Berlin, 1986, pp.~58--67. \MR{MR927758 (89b:05018)}

\bibitem{fischer}
Ilse Fischer, \emph{Another refinement of the {B}ender-{K}nuth
  (ex-){C}onjecture}, math.CO/0401235, January 2004.

\bibitem{frame}
J.~S. Frame, G.~de~B. Robinson, and R.~M. Thrall, \emph{The hook graphs of the
  symmetric groups}, Canadian J. Math. \textbf{6} (1954), 316--324.
  \MR{MR0062127 (15,931g)}

\bibitem{fulton}
William Fulton, \emph{Young tableaux}, London Mathematical Society Student
  Texts, vol.~35, Cambridge University Press, Cambridge, 1997. \MR{MR1464693
  (99f:05119)}

\bibitem{fultonharris}
William Fulton and Joe Harris, \emph{Representation theory}, Graduate Texts in
  Mathematics, vol. 129, Springer-Verlag, New York, 1991. \MR{MR1153249
  (93a:20069)}

\bibitem{gesselviennot}
Ira Gessel and G{\'e}rard Viennot, \emph{Binomial determinants, paths, and hook
  length formulae}, Adv. in Math. \textbf{58} (1985), no.~3, 300--321.
  \MR{MR815360 (87e:05008)}

\bibitem{gordon}
Basil Gordon, \emph{A proof of the {B}ender-{K}nuth conjecture}, Pacific J.
  Math. \textbf{108} (1983), no.~1, 99--113. \MR{MR709701 (85b:05019)}

\bibitem{jameskerber}
Gordon James and Adalbert Kerber, \emph{The representation theory of the
  symmetric group}, Encyclopedia of Mathematics and its Applications, vol.~16,
  Addison-Wesley Publishing Co., Reading, Mass., 1981. \MR{MR644144
  (83k:20003)}

\bibitem{knuth}
Donald~E. Knuth, \emph{Permutations, matrices, and generalized {Y}oung
  tableaux}, Pacific J. Math. \textbf{34} (1970), 709--727. \MR{MR0272654 (42
  \#7535)}

\bibitem{kostant}
Bertram Kostant, \emph{A formula for the multiplicity of a weight}, Trans.
  Amer. Math. Soc. \textbf{93} (1959), 53--73. \MR{MR0109192 (22 \#80)}

\bibitem{krattenthaler1}
C.~Krattenthaler, \emph{The major counting of nonintersecting lattice paths and
  generating functions for tableaux}, Mem. Amer. Math. Soc. \textbf{115}
  (1995), no.~552, vi+109. \MR{MR1254150 (95i:05109)}

\bibitem{krattenthaler2}
\bysame, \emph{Identities for classical group characters of nearly rectangular
  shape}, J. Algebra \textbf{209} (1998), no.~1, 1--64. \MR{MR1652177
  (2000a:05218)}

\bibitem{macmahon}
Percy~A. MacMahon, \emph{Combinatory analysis}, Two volumes (bound as one),
  Chelsea Publishing Co., New York, 1960. \MR{MR0141605 (25 \#5003)}

\bibitem{proctor}
Robert~A. Proctor, \emph{Equivalence of the combinatorial and the classical
  definitions of {S}chur functions}, J. Combin. Theory Ser. A \textbf{51}
  (1989), no.~1, 135--137. \MR{MR993658 (90b:05015)}

\bibitem{sagan}
Bruce~E. Sagan, \emph{The ubiquitous {Y}oung tableau}, Invariant theory and
  tableaux (Minneapolis, MN, 1988), IMA Vol. Math. Appl., vol.~19, Springer,
  New York, 1990, pp.~262--298. \MR{MR1035498 (90k:05017)}

\bibitem{stanley}
Richard~P. Stanley, \emph{Enumerative combinatorics. {V}ol. 2}, Cambridge
  Studies in Advanced Mathematics, vol.~62, Cambridge University Press,
  Cambridge, 1999. \MR{MR1676282 (2000k:05026)}

\end{thebibliography}
\bibliographystyle{amsplain}

\end{document}